\documentclass{cmslatex}
\usepackage[paperwidth=7in, paperheight=10in, margin=.875in]{geometry}
 \usepackage[backref,colorlinks,linkcolor=red,anchorcolor=green,citecolor=blue]{hyperref}
\usepackage{amsfonts,amssymb}
\usepackage{amsmath}
\usepackage{graphicx}
\usepackage{cite}
\usepackage{enumerate}
\renewcommand{\theequation}{\arabic{section}.\arabic{equation}}

\usepackage{amssymb, amsfonts,bbm}
\usepackage{graphicx,color}
\usepackage{diagbox}
\usepackage{algorithm}
\usepackage{algpseudocode}
\usepackage{enumitem}
\usepackage{multirow}
\usepackage{subfig}
\usepackage{epstopdf}
\usepackage{comment}
\usepackage{adjustbox}
\usepackage{todonotes}
\usepackage{url}
\usepackage{hyperref}

\newcommand{\bvec}[1]{\mathbf{#1}}

\newcommand{\vu}{\bvec{u}}

\newcommand{\vx}{\bvec{x}}

\newcommand{\vz}{\bvec{z}}

\newcommand{\Or}{\mathcal{O}}

\newcommand{\NN}{\mathbb{N}}
\newcommand{\RR}{\mathbb{R}}

\global\long\def\R{\mathbb{R}}

\numberwithin{equation}{section}
\numberwithin{figure}{section}

\newcommand{\email}[1]{\protect\href{mailto:#1}{#1}}

\allowdisplaybreaks
\begin{document}
\title{Universal approximation of symmetric and anti-symmetric functions}


\author{Jiequn Han\thanks{Department of Mathematics, Princeton University, Princeton, NJ 08544, USA; Center for Computational Mathematics, Flatiron Institute, New York, NY 10010, USA (\email{jiequnhan@gmail.com})}
\and Yingzhou Li\thanks{School of Mathematical Sciences, Fudan University,
Shanghai, China 200433 (\email{yingzhouli@fudan.edu.cn})}
\and Lin Lin\thanks{Department of Mathematics, University of California, Berkeley,  CA 94720, USA; Computational Research Division, Lawrence Berkeley National Laboratory, Berkeley, CA 94720, USA (\email{linlin@math.berkeley.edu})} 
\and Jianfeng Lu\thanks{Department of Mathematics, Department of Physics, and Department of Chemistry, Duke University, Durham, NC 27708, USA (\email{jianfeng@math.duke.edu})}
\and Jiefu Zhang\thanks{Department of Mathematics, University of California, Berkeley,  CA 94720, USA  (\email{jiefuzhang@berkeley.edu})}
\and Linfeng Zhang\thanks{Program in Applied and Computational Mathematics, Princeton University, Princeton, NJ 08544, USA (\email{linfengz@princeton.edu})}
}

\pagestyle{myheadings} \markboth{Universal Approximation of Symmetric and Anti-Symmetric Functions}{J. Han, Y. Li, L. Lin, J. Lu, J. Zhang, and L. Zhang} \maketitle

\begin{abstract}
We consider universal approximations of symmetric and anti-symmetric functions, which are important for applications in quantum physics, as well as other scientific and engineering computations. We give constructive approximations with explicit bounds on the number of parameters   with respect to the dimension and the target accuracy $\epsilon$. While the approximation still suffers from the curse of dimensionality, to the best of our knowledge, these are the first results in the literature with explicit error bounds for functions with symmetry or anti-symmetry constraints.
\end{abstract}


\section{Introduction}

We consider in this study universal approximations of symmetric and
anti-symmetric functions. A function $f:(\RR^d)^N\to \RR$ is {\bf
(totally) symmetric} if
\begin{equation} \label{eqn:symmetric}
    f(\vx_{\sigma(1)},\ldots,\vx_{\sigma(N)})=f(\vx_1,\ldots,\vx_N),
\end{equation}
for any permutation $\sigma\in S(N)$, and elements $\vx_1,\ldots,\vx_N\in
\RR^d$. Similarly $f$ is {\bf (totally) anti-symmetric} if
\begin{equation} \label{eqn:antisymmetric}
    f(\vx_{\sigma(1)},\ldots,\vx_{\sigma(N)}) = 
    (-1)^{\sigma}f(\vx_1,\ldots,\vx_N),
\end{equation}
for any permutation, 
where $(-1)^{\sigma}$ is the signature of $\sigma$.  Note that the
permutation is only applied to the outer indices $1,\ldots,N$ (also referred to as the ``particle'' indices later), 
but not the Cartesian indices $1,\ldots,d$ for each $\vx_i$. In
other words, $f$ is not totally symmetric / anti-symmetric when viewed
as a function on $\RR^{d\times N}$. More precisely, $f$ in \eqref{eqn:antisymmetric} lies in the anti-symmetric inner product of $N$ copies of $L^2(\mathbb{R}^d)$.
This is the relevant setup in
many applications in scientific and engineering computation. A totally
symmetric function is also called a permutation invariant function. A
closely related concept is the permutation equivariant mapping, which
is of the form $Y:(\RR^d)^N\to (\RR^{\tilde{d}})^N$ that satisfies
\begin{equation} \label{eqn:}
    Y_i(\vx_{\sigma(1)},\ldots,\vx_{\sigma(N)}) =
    Y_{\sigma(i)}(\vx_1,\ldots,\vx_N), \qquad i = 1, \ldots, N
\end{equation}
for any permutation $\sigma\in S(N)$, and $\vx_1,\ldots,\vx_N\in \RR^d$.
Here each component $Y_i\in\RR^{\tilde{d}}$, and $\tilde{d}$ can
be different from $d$.

Perhaps the most important example of totally symmetric and
anti-symmetric functions is the wavefunction of identical particles in
quantum mechanics. The indistinguishability of identical particles
implies that their wavefunctions should be either totally symmetric or
totally anti-symmetric upon exchanging the variables associated with
any two of the particles, corresponding to two categories of particles:
bosons and fermions.
The former can share quantum states, giving rise to, e.g., the
celebrated Bose-Einstein condensate; while the latter cannot share
quantum states as described by the famous Pauli exclusion principle.
Such exchange/permutation symmetry also arises from other applications
than identical particles in quantum mechanics, mostly in the form of
symmetric functions.  For instance, in chemistry and materials science,
the interatomic potential energy should be invariant under the permutation
of the atoms of the same chemical species. Another example is in computer
vision, where the classification of point clouds should not depend on
the ordering of points.

The dimension of symmetric and anti-symmetric functions is usually large
in practice because it is proportional to the number of considered
elements. This means that in computation the notorious difficulty of
``curse of dimensionality" is often encountered when dealing with such
functions. Recent years have witnessed compelling success of neural
networks in representing high-dimensional symmetric functions with
great accuracy and efficiency, see, e.g.,  \cite{behler2007generalized,schutt2017schnet,Zhang2018deep,Zhang2018end} for interatomic potential energy, \cite{Qi2017pointnet,Qi2017pointnet++} for
3D classification and segmentation of point sets, and \cite{germain2021deepsets,zhou2021frame} for solutions of partial differential equations. For anti-symmetric
functions, some recent work in the past few years~\cite{Choo2020fermionic,Han2019solving,Hermann2020deep,klus2021kernels,luo2019backflow,pfau2020abinitio,stokes2020phases} has shown exciting
potential of solving the many-electron Schr{\"o}dinger equation with
neural networks. Within the framework of variational Monte Carlo (VMC),
for some benchmark systems,
the anti-symmetric wavefunction parameterized by neural networks can be
on a par with the state-of-the-art wavefunctions constructed
based on chemical or physical knowledge.  

Despite the empirical success of neural networks approximating
symmetric and anti-symmetric functions, theoretical understanding of
these approximations is still limited. There are numerous results
(see e.g.~\cite{Barron1993universal,Cybenko1989approximation,Hornik1989multilayer}) concerning the universal approximation of general
continuous functions on compact domains.  Nevertheless, if the target
function is symmetric or anti-symmetric, it is much less investigated
whether one can achieve the universal approximation with a class of
functions with the same symmetry constraints.  Explicitly guaranteeing
the symmetric or anti-symmetric property of an ansatz is often mandatory.
For example, for electronic systems, if the wavefunction is not
constrained within the space of anti-symmetric functions, the resulting
variational energy could be lower than the exact ground-state energy
and would be no longer physically meaningful.  From a machine learning
perspective, symmetries can also significantly reduce the number of
effective degrees of freedom, improve the efficiency of training, and
enhance the generalizability of the model (see e.g.~\cite{han2020integrating}).
However, one needs to first make sure that the function class is still universal and sufficiently
expressive.  

The universal approximations of symmetric functions were partially
studied in~\cite{Zaheer2017deepset}. However, as will be illustrated
in Section~\ref{sec:sym1}, the proof of \cite{Zaheer2017deepset}
only holds for the case when $d=1$. Moreover, there is no error
estimation provided for the proposed approximation.  The more recent
work~\cite{Sannai2019universal} considered the universal approximation
of permutation invariant functions and equivariant mappings for $d=1$
as well. The work \cite{Sannai2019improved} gives a generalization bound that is improved by a factor of $\sqrt{N!}$ with the introduction of the quotient feature space, but the result only applies for $d=1$.
By respecting the permutation symmetry, the resulting neural
network involves much fewer parameters than the corresponding dense
neural networks.
The recent works~\cite{dusson2019atomic,bachmayr2021polynomial} on polynomial approximation of symmetric functions display a reduction of the curse of dimensionality when the target function has low-order many-body expansions or sparse polynomial approximates.
The work~\cite{zweig2020functional} studies learning symmetric functions from a different perspective by treating symmetric functions (of any size) as functions over probability measures. But such a setting is not always applicable to practical problems.
Similar to the case of symmetric functions, when $d=1$, any anti-symmetric polynomial can be factorized as the product of a Vandermonde determinant and a
totally symmetric polynomial. Such a universal representation has been
known since Cauchy \cite{Cauchy:1815}. However, to our knowledge there is no such simple factorization for $d>1$. For $d>1$, the FermiNet \cite{pfau2020abinitio} achieves good numerical results by using a linear combination of generalized Slater Determinants, and \cite[Theorem 14]{Hutter2020representing} shows the universality of the FermiNet Ansatz again based on a global lexicographic ordering for $d>1$.  

In this paper we aim to study the universal approximation of general
symmetric and anti-symmetric functions for any $d\ge 1$. We now summarize the main results of the paper: First, for the symmetric function, we give
two different proofs of the universality of the ansatz proposed in
\cite{Zaheer2017deepset}, both with explicit error bounds.
The first proof is based on the Ryser formula~\cite{Ryser:1963} for
permanents, and the second is based on the partition of the state space,
as elaborated in Section~\ref{sec:sym1} and Section~\ref{sec:sym2},
respectively.  Moreover, we also show in Section~\ref{sec:antisym}
that, for the general anti-symmetric function with elements in any
dimension, a simple ansatz combining Vandermonde determinants and the
ansatz for symmetric functions is universal, with similar explicit
error bounds as in the symmetric case.  For readers' convenience, we
summarize below in two theorems the ansatz we proved with the universal
approximation for symmetric and anti-symmetric functions. 
The approximation rate only relies on a weak condition that the gradient is uniformly bounded.
We conclude in Section~\ref{sec:practical}
with some practical considerations and future directions for further
investigation. The proofs of Theorem \ref{thm:sym} and \ref{thm:antisym} are given in Section \ref{sec:sym2} and \ref{sec:antisym}, respectively.


\begin{theorem}[Approximation to  symmetric functions]
Let $f: \Omega^{N}\to \R$ be a continuously differentiable, totally symmetric function, where $\Omega$  is a compact subset of $\RR^d$. Let  $0<\epsilon<\|\nabla f\|_2 \sqrt{Nd} N^{-\frac{1}{d}}$, here $\|\nabla f\|_2 = \max_X \|\nabla f(X)\|_2$. Then there exist $g:\RR^d\to\RR^M$, $\phi:\RR^M\to\RR$, such that for any $X=(\vx_1,\ldots,\vx_N)\in \Omega^N$,
\begin{displaymath}
\left|f(X)-\phi\left(\sum_{j=1}^N g(\vx_j)\right)\right| \le \epsilon,
\end{displaymath}
where $M$, the number of feature  variables, satisfies the bound
\begin{equation}
M \leq 2^N(\|\nabla f\|_2^2Nd)^{Nd/2} / (\epsilon^{Nd}N!).
\label{eqn:bound_term}
\end{equation}
\label{thm:sym}
\end{theorem}

\begin{theorem}[Approximation to anti-symmetric functions]
 Let $f: \Omega^{N}\to \R$ be a continuously differentiable, totally anti-symmetric function, where $\Omega$  is a compact subset of $\RR^d$. Then  there exist $K$ permutation equivariant mappings  $Y^{k}:(\RR^d)^{N}\to\RR^{N}$, and permutation invariant functions $U^k:(\RR^d)^N\to\RR$, $k=1,\ldots, K$, such that for any $X=(\vx_1,\ldots,\vx_N)\in \Omega^N$,
\begin{displaymath}
\left|f(X) - \sum_{k=1}^K U^k(X)\prod_{i<j}(y^k_i(X) - y^k_j(X))\right|<\epsilon,
\end{displaymath}
where $K$ satisfies the bound
\begin{displaymath}
K \leq (\|\nabla f\|_2^2Nd)^{Nd/2} / (\epsilon^{Nd}N!).
\end{displaymath}
For each  $U^k$, there exists $g^k:\RR^d\to\RR^m$, $\phi^k:\RR^m\to\RR$ with  $m\le 2^N$, such that for any $X=(\vx_1,\ldots,\vx_N)\in \Omega^N$,
\[
U^k(X)=\phi^k\left(\sum_{j=1}^N g^k(\vx_j)\right).
\]
\label{thm:antisym}
\end{theorem}

We remark that for the bounds in both theorems above, $\|\nabla f\|_2$ may also scale with respect to $N$. For the examples listed below, we assume the domain for the function to be a hypercube $\Omega_N=[0,1]^N$. In the symmetric case, if $f$ is in the simple form $f(X)=\sum_{i=1}^N g(\vx_i)$, $\|\nabla f\|_2$ will scale as $\Or(\sqrt{N})$. In the antisymmetic case, if $f$ is a Slater determinant, the integral on the global domain $\int_{\Omega^N} \|\nabla f(X)\|_2^2 dX$ (which can be interpreted as the kinetic energy of the quantum system) is expected to scale as $\Or(N)$. Hence $\|\nabla f(X)\|_2$ scales as $\Or(\sqrt{N})$ at a given $X$ on average. Both the simple examples above have $\Or(\sqrt{N})$ scaling, but the exact scaling will depend on the choice of function $f$. We also remark that Theorem \ref{thm:antisym} only relies on a local lexicographic ordering procedure, instead of a global sorting operation.

\section{Totally symmetric functions}\label{sec:sym1}

Let $X=(\vx_1,\ldots,\vx_N)$, with $\vx_i\in \Omega:=[0,1]^d$. Consider a totally symmetric function $f(X)$. It is proved in \cite{Zaheer2017deepset} that  when $d=1$ (therefore $\Omega=[0,1]$), the following universal approximation representation holds
\begin{equation}
f(X)=\phi\left(\sum_{j=1}^N g(x_j)\right)
\label{eqn:symmetry_representation}
\end{equation}
for continuous  functions $g:\RR\to\RR^M$, and $\phi:\RR^M\to\RR$. For completeness we briefly recall the proof.

Let $\mathcal{X}=\{(x_1,\ldots,x_N)\in \Omega^N = \RR^N: x_1 \leq x_2
\leq \cdots \leq x_N\}$. Define the mapping $E:\mathcal{X}\to
\mathcal{Z}\subset\RR^{N+1}$, with each component function defined as
$$
z_{q} =E_{q}(X) := \sum_{n=1}^N (x_n)^{q}, \quad q=0,1,\ldots,N.
$$
The mapping $E$ is a homeomorphism between $\mathcal{X}$ and its image in $\RR^{N+1}$~\cite{Zaheer2017deepset}. Hence, if we let $g:\RR \to \RR^{N+1} , x \mapsto [1,x,x^2,\ldots,x^N]$ and $\phi: \RR^{N+1} \to \RR, Z \mapsto f(E^{-1}(Z))$, then we have
$$ 
f(X)=\phi\left(\sum_{j=1}^N g(x_j)\right).
$$
Here the number of feature variables is $M=N+1$ by construction.
The main difficulty associated with this construction is that the mapping $E^{-1}$ can be arbitrarily complex to be approximated in practice. In fact the construction is similar in flavor to the Kolmogorov-Arnold representation theorem \cite{Kolmogorov1957}, which provides a universal representation for multivariable continuous functions, but without any \textit{a priori} guarantee of the accuracy with respect to the number of parameters.

In order to generalize to the case $d>1$, the proof of \cite[Theorem 9]{Zaheer2017deepset} in fact suggested an alternative  proof for the case $d=1$ as follows. Using the Stone-Weierstrass theorem, a totally symmetric function can be approximated by a polynomial. After symmetrization, this polynomial becomes  a totally symmetric polynomial. By the fundamental theorem of symmetric polynomials \cite{Macdonald1998}, any symmetric polynomial can be represented by a polynomial of elementary symmetric polynomials. In other words, for any symmetric polynomial $P$, we have 
$$
P(X)=Q(e_1(X),\ldots,e_N(X)),
$$ 
where $Q$ is some polynomial, and the elementary polynomials $e_k$ are defined as $$e_k(X) = \sum_{1 \leq j_1<j_2<\cdots<j_k\leq N} x_{j_1}x_{j_2}\cdots x_{j_k}.$$
Using the Newton-Girard formula, an elementary symmetric polynomial can be represented  with power sums by
\begin{equation}
e_k(X) = \frac{1}{k!} \begin{vmatrix}
E_1(X) & 1 & 0 & 0 & \cdots & 0\\
E_2(X) & E_1(X) & 2 & 0 & \cdots & 0\\
E_3(X) & E_2(X) & E_1(X) & 3 & \cdots & 0\\
\vdots & \vdots & \vdots & \vdots & \ddots & \vdots\\
E_{k-1}(X) & E_{k-2}(X) & E_{k-3}(X) & E_{k-4}(X) & \cdots & k-1\\
E_k(X) & E_{k-1}(X) & E_{k-2}(X) & E_{k-3}(X) & \cdots & E_1(X)\\
\end{vmatrix}.
\label{eqn:newton_girard}
\end{equation}
Now let $g$ be the same function defined in the previous proof, and define $\phi$ in terms of $Q$ and the determinant computation. We can obtain a polynomial approximation in the form of
$$ 
f(x) = P(X) +\epsilon = \phi\left(\sum_{j=1}^N g(x_j)\right) + \epsilon. 
$$
Here the error $\epsilon$ is due to the Stone-Weierstrass approximation. Letting $\epsilon\to 0$ we obtain the desired representation. 

However, it is in fact not straightforward to extend the two proofs above to the case $d>1$. For the first proof, we can not define an ordered set $\mathcal{X}$ when each $\vx_i\in\RR^d$ to define the homeomorphism $E$. For the second proof, in the case $d>1$ a monomial (before symmetrization) takes the form
\begin{displaymath}
\prod_{i=1}^N \prod_{\alpha=1}^d \vx_{i,\alpha}^{\gamma_{i,\alpha}}, \quad \gamma_{i,\alpha}\in \NN.
\end{displaymath}
Note that $f(X)$ is only symmetric with respect to the particle index $i$, but not the component index $\alpha$. Hence the symmetrized monomial is \textit{not} a totally symmetric function with respect to all variables. Therefore the 
fundamental theorem of symmetric polynomials does not apply.
 
Below we prove that the representation \eqref{eqn:symmetry_representation} indeed holds for any $d\ge 1$, and therefore we complete the proof 
of \cite{Zaheer2017deepset}.
For technical reasons to be illustrated below, and without loss of generality, we shift the domain $\Omega$ and assume $\vx_i\in \Omega:=[1,2]^d$. Following the Stone-Weierstrass theorem and after symmetrization, $f(X)$ can be approximated by a symmetric polynomial. Every symmetric polynomial can be written as the linear combination of symmetrized monomials of the form
$$
\sum_{\sigma\in S(N)} \prod_{i=1}^N \prod_{\alpha=1}^d \vx_{\sigma(i),\alpha}^{\gamma_{i,\alpha}}:=\sum_{\sigma\in S(N)} \prod_{i=1}^N f_i(\vx_{\sigma(i)}):=\text{perm} ([f_i(\vx_j)]).
$$
 Here $f_i(\vx_j)=\prod_{\alpha=1}^d \vx_{j,\alpha}^{\gamma_{i,\alpha}}$, $[f_i(\vx_j)]$ denotes an $N\times N$ matrix whole entry in the $i$-th row and $j$-column is $f_i(\vx_j)$, and $\text{perm}(A)$ stands for the permanent of square matrix $A$.  

Following the Ryser formula \cite{Ryser:1963} for representing a permanent (noting that permanent is invariant under transposition), we have\begin{align*}
\text{perm}([f_i(\vx_j)])=&(-1)^{N} \sum_{S \subseteq\{1, \ldots, N\}}(-1)^{|S|} \prod_{i=1}^{N} \sum_{j \in S} f_j(\vx_i)\\
=& (-1)^{N} \sum_{S \subseteq\{1, \ldots, N\}}(-1)^{|S|} e^{\sum_{i=1}^{N} \log\left(\sum_{j \in S} f_j(\vx_i)\right)}.
\end{align*}
Here we have used that  $f_j(\vx)>0$ for all $\vx\in\Omega$. Now we write down the approximation using a symmetric polynomial, which is a linear combination of $L$ symmetrized monomials
\begin{align*}
f(X)-\epsilon = P(X) =& \sum_{l=1} ^{L}c^{(l)} \text{perm} ([f^{(l)}_i(\vx_j)]) \\
=& (-1)^N \sum_{l=1}^L c^{(l)}\sum_{S \subseteq\{1, \ldots, N\}}(-1)^{|S|} e^{\sum_{i=1}^{N} \log\left(\sum_{j \in S} f^{(l)}_j(\vx_i)\right)}
\end{align*}
Define $g:\RR^d \to \RR^{L2^N}$ with each component function
$$
g^{(l)}_S(\vx) = \log\left(\sum_{j \in S} f^{(l)}_j(\vx)\right).
$$
Then we define $\phi: \RR^{L2^N} \to \RR$ given by
$$
\phi(Y) = (-1)^N \sum_{l=1} ^{L}c^{(l)}\sum_{S \subseteq\{1, \ldots, N\}}(-1)^{|S|} e^{Y^{(l)}_S},
$$
where $Y^{(l)}_S$ is the $S$-th component of $g^{(l)}$. 
We now have an approximation of the target totally symmetric function in the desired form
$$
f(X)=\phi\left(\sum_{j=1}^N g(\vx_j)\right)+\epsilon,
$$
and we finish the proof. Here the number of feature variables is $M = L2^N$, where $L$ is the number of symmetrized monomials used in the approximation. 
\section{Totally symmetric function, revisited}\label{sec:sym2}

In this section, we  prove Theorem \ref{thm:sym} for any $d\ge 1$. In particular, our proof is more explicit and does not rely on the Stone-Weierstrass theorem. 
The main idea is to partition the space into a lattice and use piecewise-constant functions to approximate the target permutation invariant function.

Again without loss of generality we assume $\vx_i\in[0,1]^d:=\Omega$. We
then partition the domain $\Omega$ into a lattice $\mathbb{L}$
with grid size $\delta$ along each direction. Due to symmetry, we
can assign a lexicographical order $\preceq$ to all lattice points
$\vz_i\in\mathbb{L}$. That is, $\vz_1 \preceq \vz_2$ if $\vz_{1,\alpha} <
\vz_{2,\alpha}$ for the first $\alpha$ where $\vz_{1,i}$ and $\vz_{2,i}$
differs, $\alpha = 1,2,\cdots,d$.  We define the tensor product of the $N$ copies of the lattice $\mathbb{L}$ as $\mathbb{L}^N$, and a wedge of  $\mathbb{L}^N$ is defined accordingly as
$$
\bigwedge\nolimits^{\!N}\mathbb{L}:=\{Z=(\vz_1,\ldots,\vz_N)|\vz_1\preceq\vz_2\cdots\preceq\vz_N\}.
$$
For each $Z\in \bigwedge^N \mathbb{L}$, a corresponding union of boxes in $\Omega^N$ can be written as 
$$
B^{Z;\delta}=\bigcup_{\sigma\in S(N)}\{X~\vert~\vx_i=\vz_{\sigma(i)}+\delta \vu_i, \quad \vu_i\in[0,1]^d\}.
$$
By construction, the piecewise-constant approximation to the target permutation invariant function is then 
$$
f(X) =  \sum_{Z \in\bigwedge^N \mathbb{L}}  f(Z) \mathbbm{1}_{B^{Z;\delta}}(X)+\epsilon.
$$
We assume that the derivative $\nabla_X f(X)$ is uniformly bounded for $X\in \Omega^N$ and denote the bound in 2-norm by $\|\nabla f\|_2$. The maximal distance between two points in a box is bounded by the length of the longest diagonal, which is $\delta\sqrt{Nd}$. Hence the approximation error satisfies $|\epsilon| \le \|\nabla f\|_2\delta\sqrt{Nd}$, which is obtained by applying the mean value theorem and the Cauchy-Schwarz inequality. Note that the indicator function $\mathbbm{1}_{B^{Z;\delta}}(X)$ is permutation invariant and can be rewritten as
\begin{align*}
\mathbbm{1}_{B^{Z;\delta}}(X) =& \frac{1}{C_Z}\sum_{\sigma\in S(N)} \mathbbm{1}_{\{X|\vx_i=\vz_{\sigma(i)}+\delta \vu_i, \vu_i\in[0,1]^d,1\le i\le N\}}(X)\\
=& \frac{1}{C_Z}\sum_{\sigma\in S(N)} \prod_{i=1}^N \mathbbm{1}_{\{\vx_i=\vz_{\sigma(i)}+\delta \vu_i, \vu_i\in[0,1]^d\}}(\vx_i)\\
=& \frac{1}{C_Z}\sum_{\sigma\in S(N)} \prod_{i=1}^N \mathbbm{1}_{\{\vx_{\sigma(i)}=\vz_i+\delta \vu_i, \vu_i\in[0,1]^d\}}(\vx_{\sigma(i)}) = \frac{1}{C_Z} \text{perm} ([f^Z_i(\vx_j)]),
\end{align*}
where $f^Z_i(\vx)=\mathbbm{1}_{\{\vx|\vx=\vz_{i}+\delta\vu, \vu\in[0,1]^d\}}(\vx)$. The constant $C_Z$ takes care of repetition that can happen depending on $Z$. When all elements in $Z$ are distinct, the box $X$ lives in only corresponds to one permutation, so $C_Z=1$ in this case. If say $\vz_1=\vz_2$ and all other elements distinct, then the box $X$ lives in has two corresponding permutations that differ by a swapping of the first two elements. In this case, $C_Z=2$ will account for the arising repetition.
Next we apply the Ryser formula to the permanent,
$$
\mathbbm{1}_{B^{Z;\delta}}(X)=\frac{1}{C_Z}\text{perm} ([f^Z_i(\vx_j)])= \frac{(-1)^N}{C_Z} \sum_{S \subseteq\{1, \ldots, N\}}(-1)^{|S|} e^{\sum_{i=1}^{N} \log\left(\sum_{j \in S} f^Z_j(\vx_i)\right)}.
$$
We can now define $g: \RR^d \to \RR^{|S(N)|\times|\bigwedge^N \mathbb{L}|}$ where each component function is given by $$g_S^Z(\vx)=\log\left(\sum_{j \in S} f^Z_j(\vx)\right)$$
and we define $\phi: \RR^{|S(N)|\times|\bigwedge^N \mathbb{L}|} \to \RR$ as
\begin{equation}
\phi(Y) = \sum_{Z \in\bigwedge^N \mathbb{L}}  \frac{(-1)^N f(Z)}{C_Z}  \sum_{S \subseteq\{1, \ldots, N\}}(-1)^{|S|} e^{Y^Z_S}.
\label{eqn:sym_approx_phi}
\end{equation}
Since $f^Z_j$'s are indicator functions we naturally have $\sum_{j\in S} f^Z_j(\vx)\geq 0$. In the case when $\sum_{j\in S} f^Z_j(\vx)=0$, $g_S^Z(\vx)=-\infty$. In this case, $e^{g_S^Z(\vx)}=0$, and therefore its contribution to  $\phi$ vanishes as desired. In summary, we arrive at the universal approximation
\begin{equation}
f(X) = \phi\left(\sum_{j=1}^N g(\vx_j)\right)+\epsilon.
\label{eqn:sym_approx_partition}
\end{equation}

Due to the explicit tabulation strategy, the number of terms needed in
the approximation \eqref{eqn:sym_approx_partition} can be counted as
follows. The number of points in $\bigwedge\nolimits^{\!N}\mathbb{L}$
is $\Or((\delta^{-Nd})/N!)$, where $N!$ comes from the lexicographic
ordering. Note that formally as $N\to\infty$, $(\delta^{-Nd})/N!\sim
(\delta^{-d}/N)^{N}$ can vanish for fixed $\delta$. However, this
means that the number of elements has exceeded the number of grid
points in $\mathbb{L}$ and is unreasonable. So we should at least have
$\delta\lesssim N^{-\frac{1}{d}}$. 
In order to obtain an $\epsilon$-close
approximation of $f(X)$, we require our error bound $\|\nabla f\|_2\delta \sqrt{Nd} \sim \epsilon$. 
When $\epsilon\le \|\nabla f\|_2\sqrt{Nd} N^{-\frac{1}{d}}$, we can choose $\delta = \frac{\epsilon}{\|\nabla f\|_2} (Nd)^{-\frac{1}{2}}$ so that both $\delta\lesssim N^{-\frac{1}{d}}$ and $\|\nabla f\|_2 \delta \sqrt{Nd} \sim \epsilon$ are fulfilled, and the
number of points in $\bigwedge\nolimits^{\!N}\mathbb{L}$ becomes
$(\|\nabla f\|_2^2Nd)^{Nd/2} / (\epsilon^{Nd}N!)$ because $\delta^{-Nd} = (\|\nabla f\|_2^2Nd)^{Nd/2} / \epsilon^{Nd}$.
For each $Z$, the
number of terms to be summed over in Eq. \eqref{eqn:sym_approx_phi} is
$|S(N)| = 2^N$. Therefore in order to obtain an $\epsilon$-approximation,
the number of feature variables is given by Eq. \eqref{eqn:bound_term}.
This proves Theorem~\ref{thm:sym}.

This is  of course a very pessimistic bound, and we will discuss on
the practical implications for designing neural network architectures
in Section~\ref{sec:practical}.  We remark that one may expect
that following the same tabulation strategy, we may also provide a
quantitative bound for $\phi$ constructed by the homeomorphism mapping
$E^{-1}$ as discussed in Section~\ref{sec:sym1}. However, the difference
is that our bound only relies on the smoothness of the original
function $f$ and hence the bound for $\delta$. On the other hand,
the mapping $E^{-1}$ and hence $\phi$ can be arbitrarily pathological,
and therefore it is not even clear how to obtain a double-exponential
type of bound as discussed above. 
We also remark that if the 
indicator functions $\mathbbm{1}_{B^{Z;\delta}}(X)$ and 
$\mathbbm{1}_{\{\vx|\vx=\vz_{i}+\delta\vu, \vu\in[0,1]^d\}}(\vx)$ in the proof are replaced by proper smooth cutoff functions with respect to corresponding domains, the ansatz in Eq. \eqref{eqn:sym_approx_partition} 
can be continuous to accommodate the applications that require continuity.
For example, an interatomic potential energy
should be continuous to guarantee the total energy is conserved during
molecular dynamics simulations.

\section{Totally anti-symmetric functions}\label{sec:antisym}

Now we consider an anti-symmetric function. Similar to the symmetric
case, when $d=1$ and $f(X)$ is a polynomial of $X$ and anti-symmetric,
it is known that
$$
f(x_1,\ldots,x_N) = U(X) \prod_{i<j} (x_i-x_j),
$$
where $U(X)$ is a symmetric polynomial and the second term is
a Vandermonde determinant. This was first proved by Cauchy
\cite{Cauchy:1815}, who of course also first introduced the concept of
determinant in its modern sense.

Our aim is to generalize to $d > 1$. Without loss of generality we again assume $\vx_i\in[0,1]^d:=\Omega$. 
The construction of the ansatz is parallel to the totally symmetric case. 
Recall the lattice $\mathbbm{L}$, the wedge $\bigwedge^N \mathbbm{L}$, and the corresponding union of boxes $B^{Z, \delta}$.
In the totally symmetric case, we make a piece-wise constant approximation over the union of boxes. For the anti-symmetric case, we would need to insert an anti-symmetric factor (w.r.t. $X$): 
\begin{equation}
  f(X) = \sum_{Z \in \bigwedge^N\mathbb{L}} f(Z) \frac{\psi^Z(X)}{\psi^Z(Z)} \mathbbm{1}_{B^{Z;\delta}}(X)+\epsilon, 
\end{equation}
where $\psi^Z(\cdot)$ is a totally anti-symmetric function which might be chosen depending on $Z$. Here the approximation error, similar as before, satisfies $|\epsilon| < \|\nabla f\|_2\delta\sqrt{Nd}$. Motivated by the one-dimension result, we focus on constructions of $\psi^Z$ given by the Vandermonde determinant. Given a permutation equivariant map $X \mapsto Y^Z: (\RR^{d})^N \to \RR^N$, we consider 
\begin{equation}
  \psi^Z(X) = \prod_{i<j}(y_i(X) - y_j(X)). 
\end{equation}
It is clear that due to the permutation equivariance of $Y^Z$, the $\psi^Z$ defined is anti-symmetric. 
Observe that due to anti-symmetry $f(Z) = 0$ whenever $\bvec{z}_i = \bvec{z}_j$ for some $i \neq j \in [n]$. In particular, this means that we only need to consider $\psi^Z$ for those $Z$ that all $\bvec{z}_i$'s are different. 
As $\bvec{z}_i$'s are distinct, for $X \in B^{Z; \delta}$, there exists a unique $\sigma \in S(N)$ such that 
\begin{equation}
  X \in \bigl\{ X' \mid \bvec{x}_i' = \bvec{z}_{\sigma(i)} + \delta \bvec{u}_i, \quad \bvec{u}_i \in [0, 1]^d \bigr\}.
\end{equation}
Denote this unique permutation as $\sigma_{Z, X}$, and we take the permutation equivariant map $Y^Z$ such that 
\begin{equation}
\label{eq:construction1Ymap}
    y_i^Z(X) = \sigma_{Z, X}(i).
\end{equation}
Since $\sigma_{Z, X}$ gives the sorting of $X$ according to $Z$, it is easy to see that the above $Y^Z$ is equivariant, and 
\begin{equation}
    \psi^Z(X) = \prod_{i<j} (y_i(X) - y_j(X)) =  (-1)^{\sigma_{Z, X}} \prod_{i<j} (i - j).  
\end{equation}
Note that $\psi^Z(Z) = \prod_{i<j} (i - j)$, so we arrive at \begin{equation}
  f(X) = \sum_{Z \in \bigwedge^N \mathbb{L}} U^Z(X) \prod_{i<j} (y_i(X) - y_j(X))  +\epsilon.
\end{equation}
with 
\begin{equation*}
  U^Z(X) = \frac{f(Z)}{\prod_{i<j} (i -j)} \mathbbm{1}_{B^{Z;\delta}}(X).
\end{equation*}
In this construction, $U^Z(X)$ is a scaled version of the characteristic function, and can thus be treated similarly as in the totally symmetric case based on Ryser's formula. 
As the construction depends on the same tabular strategy as the totally symmetric case, the number of terms involved in the sum over $\bigwedge^N \mathbbm{L}$ is the same too, which we will not repeat. The total number of feature variables is the same as that in Eq.  \eqref{eqn:bound_term} due to the use of the indicator function. This proves Theorem~\ref{thm:antisym}. 

We remark that the index-based permutation equivariant map $Y^Z$ in \eqref{eq:construction1Ymap} leads to discontinuity at certain points. 
We may modify this permutation equivariant map in order to make the whole approximation continuous so that it is more suitable for scientific applications such as representing the many-body wavefunctions in VMC.
However, such a continuous approximation in the current tabular strategy would imply that the anti-symmetric function $\psi^Z(\cdot)$ has the nodes close to those of $f$, which seems to be a very strong assumption. We leave the error bounds of a continuous approximation to an anti-symmetric function to future works.

\section{Practical considerations and discussion}\label{sec:practical}

In this paper we study the universal approximation for symmetric and anti-symmetric functions. Following the line of learning theory, there are many questions open. For instance, the impact of symmetry on the generalization error remains unclear. This requires in-depth understanding of the suitable function class for symmetric and anti-symmetric functions, such as some adapted Barron space~\cite{Ma2019barron}. 
Note that a recent work \cite{Sannai2019improved} investigates the approximation and generalization bound of permutation invariant deep neural networks in the general case $d \geq 1$,
however with two limitations. The first is a rather strong assumption
that the target function is Lipschitz with respect to the $\ell_\infty$
norm (but not the usual Euclidean norm). This can be a severe limitation as the dimension (both $N$ and $d$) increases.
Indeed, under the same Lipschitz assumption of \cite{Sannai2019improved}, the number of feature variables $M$ in our Theorem~\ref{thm:sym} can be improved accordingly to $\Or\left(2^N/ (\epsilon^{Nd}N!)\right)$.
The second limitation is that
the proposed ansatz in \cite{Sannai2019improved} introduces \textit{sorting layers} to represent the sorting procedure
at the first step. The sorting procedure will bring discontinuity,
which leads to serious problems in some scientific applications, such as interatomic potential energy in molecular dynamics simulations.


For anti-symmetric function, 
the ansatz suitable for the practical computation is of interest since one main motivation for studying anti-symmetric function is to integrate neural network-based wavefunction into VMC. 
The Vandermonde determinant considered here is proved to provide a simple but universal ansatz. 
Its universality suggests us to consider the following trial wavefunction in VMC:
\begin{equation*}
  f(X) = \sum_{k=1}^K U^k(X)\prod_{i<j}(y^k_i(X) - y^k_j(X)),
\end{equation*}
where $U^k(X)$ is symmetric and $(y^k_1(X),\ldots,y^k_N(X))$ is a  permutation equivariant map. The ansatz for each $U^k$ and $Y^k$ 
 can be still quite flexible. 
Another more general yet more complicated ansatz is based on replacing the above Vandermonde determinant 
with the Slater determinant (see e.g. \cite{pfau2020abinitio,Hermann2020deep})
\begin{equation*}
\det [Y^k(X)] = 
\begin{vmatrix}
y^k_1(X) & \ldots & y^k_1(X) \\
\vdots & & \vdots \\
y^k_N(X) & \ldots & y^k_N(X)
\end{vmatrix}.
\end{equation*}
The Slater determinant is widely used in quantum chemistry. It is known to be universal under the complete basis set 
and, indeed, the basis set derived from the Hartree-Fock approximation often provides a fairly good starting point for most of the modern quantum chemistry methods.
However, the complexity of computing a Slater determinant is $\Or(N^3)$, while it is only $\Or(N^2)$ for a Vandermonde determinant.
This may become a more severe issue when one calculates the local energy, which involves the evaluation of the Laplacian of the trial wavefunction.
Therefore, it remains interesting to spend more effort on the Vandermonde ansatz and its variants, hoping to find good ones that strike a good balance between accuracy and efficiency. It would also be interesting to learn from the second quantized representation of quantum systems, which lifts the symmetry requirement 
of functions to that of linear operators, and leads to powerful representations such as matrix product states.

\par{\bf Acknowledgements.}
We thank the hospitality of the American Institute of Mathematics (AIM) for the workshop ``Deep learning and partial differential equation'' in October 2019, which led to this collaborative effort. The work of LL and JZ was supported in part by  the Department of Energy under Grant No. DE-SC0017867, and No. DE-AC02-05CH11231.
The work of YL and JL was also supported in part by the National Science Foundation via grants DMS-1454939 and ACI-1450280.






\bibliographystyle{spmpsci}      

\end{document}